\documentclass[hyper,11pt]{JHEP3}

\usepackage{xy,amssymb,amsmath,diagrams,times,theorem,longtable}

\xyoption{poly}
\xyoption{arc}
\xyoption{knot}
\xyoption{curve}
\xyoption{arrow}
\xyoption{all}

\newcommand{\C}{\mathbb{C}}

\newcommand{\N}{\mathbb{N}}
\newcommand{\FF}{\mathbb{F}}
\newcommand{\Z}{\mathbb{Z}}

\newcommand{\vtx}[1]{*+[o][F-]{\scriptscriptstyle #1}}

\newcommand{\wis}[1]{{\text{\em \usefont{OT1}{cmtt}{m}{n} #1}}}

\newtheorem{theorem}{Theorem}

\newtheorem{lemma}{Lemma}

{\theorembodyfont{\rmfamily}

 }

\preprint{UIA preprint 2003-03}

\title{Canonical systems and non-commutative geometry}

\author{Lieven Le Bruyn \\
Departement Wiskunde en Informatica, Universiteit Antwerpen  \\
B-2020 Antwerp (Belgium) \\
E-mail :  \email{lieven.lebruyn@ua.ac.be}}

\author{Markus Reineke \\
Bergische Universit\"at Wuppertal, Gau\ss str. 20 \\
D-42097 Wuppertal \\
E-mail : \email{reineke@math.uni-wuppertal.de}}

\abstract{ Inspired by ideas from non-commutative geometry, unions of  moduli spaces of linear control systems are identified as open subsets of infinite Grassmannians.}

\begin{document}

\begin{center} {\it For Michiel Hazewinkel on his 60th birthday.}
\end{center}

\section{Introduction}

A {\em linear control system} $\Sigma$ of type $(m,n,p)\in\N^3$ is determined by the system of linear differential equations
\[
\begin{cases}
\frac{dx}{dt} &= A x + B u \\
y &= C x,
\end{cases}
\]
where $u(t) \in \C^m$ is the {\em input or control} at time $t$, $x(t) \in \C^n$ is the {\em state} of the system and $y(t) \in \C^p$ is its {\em output}. That is, $\Sigma$ is described by a triple of matrices
\[
\Sigma = (A,B,C) \in M_n(\C) \times M_{n \times m}(\C) \times M_{p \times n}(\C) = V_{m,n,p} \]
and is said to be equivalent to a system $\Sigma' = (A',B',C') \in V_{m,n,p}$ if and only if there is a basechange matrix $g \in GL_n=GL_n(\C)$ in the state-space such that 
\[
\Sigma \sim \Sigma' \quad \Leftrightarrow \quad A' = gAg^{-1},\quad B'=gB \quad \text{and} \quad C' = C g^{-1}. \]
A system $\Sigma = (A,B,C) \in V_{m,n,p}$ is said to be {\em completely controllable} (resp. {\em completely observable}) if and only if the matrix
\[
c(\Sigma) = \begin{bmatrix} B & AB & A^2B & \hdots & A^{n-1}B \end{bmatrix} \quad ( \text{resp.} \quad
o(\Sigma) = \begin{bmatrix} C \\ CA \\ CA^2 \\ \vdots \\ CA^{n-1} \end{bmatrix}~) \]
is of maximal rank. These conditions define $GL_n$-open subsets $V_{m,n,p}^{cc}$, resp. $V_{m,n,p}^{co}$, consisting of systems with trivial $GL_n$-stabilizer, whence we have corresponding orbit spaces
\[
\wis{sys}_{m,n,p}^{cc} = V_{m,n,p}^{cc}/GL_n \qquad \text{and} \qquad \wis{sys}_{m,n,p}^{co} = V_{m,n,p}^{co}/GL_n, \]
which are known to be smooth quasi-projective varieties of dimension $(m+p)n$, see for example \cite[Part IV]{Tannenbaum}. A system $\Sigma = (A,B,C) \in V_{m,n,p}$ is said to be {\em canonical} if it is both completely controllable and completely observable. The corresponding moduli space
\[
\wis{sys}_{m,n,p}^{c} = (V_{m,n,p}^{cc} \cap V_{m,n,p}^{co})/GL_n \]
classifies canonical systems having the same input-output behavior, that is, such that all the $p \times m$ matrices $C A^i B$ for $i \in \N$ are equal \cite[Part VI - VII]{Tannenbaum}. Conversely, if $F = \{ F_j~:~j \in \N_+ \}$ is a sequence of $p \times m$ matrices such that the corresponding {\em Hankel matrices}
\[
H_{ij}(F) = \begin{bmatrix}
F_1 & F_2 & \hdots & F_j \\
F_2 & F_3 & \hdots & F_{j+1} \\
\vdots & \vdots & & \vdots \\
F_i & F_{i+1} & \hdots & F_{i+j-1}
\end{bmatrix}
\]
are such that there exist integers $r$ and $s$ such that $rk~H_{rs}(F) = rk~H_{r+1,s+j}(F)$ for all $j \in \N_+$, then $F$ is {\em realizable} by a canonical system $\Sigma = (A,B,C) \in V_{m,n,p}^c$ (for some $n$ which is equal to $rk~H_{rs}(F)$), that is,
\[
F_j = CA^{j-1}B \qquad \text{for all $j \in \N_+$}, \]
see for example \cite[Part VI - VII]{Tannenbaum} for connections between this realization problem and classical problems in analysis. These problems would be facilitated if there was an infinite dimensional manifold $X$ together with a natural stratification
\[
X = \bigsqcup_n \wis{sys}_{m,n,p}^c \]
by the moduli spaces of canonical systems (for fixed $m$ and $p$ and varying $n$).

Non-commutative geometry, as outlined by M. Kontsevich in \cite{Kontsevich}, offers a possibility to glue together closely related moduli spaces into an infinite dimensional variety controlled by a non-commutative algebra. The individual moduli spaces are then recovered as moduli spaces of simple representations (of specific dimension vectors) of the non-commutative algebra. An illustrative example is contained in the recent work by G. Wilson and Yu. Berest \cite{Wilson} \cite{BerestWilson} relating Calogero-Moser spaces to the adelic Grassmannian (see also \cite{LeBruynBocklandt} and \cite{Ginzburg} for the connection with non-commutative geometry). The main aim of the present paper is to offer another (and more elementary) example :

\begin{theorem} The equivalence classes of canonical systems with fixed input- and output-dimensions $m$ and $p$ form a specific open submanifold
\[
\bigsqcup_n \wis{sys}_{m,n,p}^c \rInto \wis{Gras}_{m+p}(\infty) \]
of the infinite Grassmannian of $m+p$-dimensional subspaces.
\end{theorem}

This paper is organized as follows. In section two we show that Kontsevich's approach is applicable to moduli spaces of canonical systems by proving that there is a natural one-to-one correspondence between equivalence classes of canonical systems with $n$-dimensional state space and isomorphism classes of simple representations of dimension vector $(1,n)$ of the formally smooth path algebra of the quiver
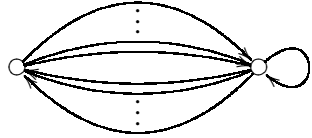
\[
\xymatrix@=3cm{
\vtx{} \ar@/^/[r] \ar@/^2ex/[r]^{\vdots} \ar@/^5ex/[r] & \vtx{} \ar@(ur,dr) \ar@/^/[l] \ar@/^2ex/[l] \ar@/^5ex/[l]_{\vdots} }
\]
with $m$ arrows pointing right and $p$ arrows pointing left. This observation gives a short proof of the following result, due to M. Hazewinkel (\cite[thm. VI.2.5]{Tannenbaum} or \cite[(2.5.7)]{Hazewinkel2}):

\begin{theorem}[Hazewinkel] The moduli space  $\wis{sys}_{m,n,p}^c$ of canonical systems is a quasi-affine variety.
\end{theorem}

In section 3 we prove that the moduli spaces $\wis{sys}_{m,n,p}^{cc}$ (resp. $\wis{sys}_{m,n,p}^{co}$) of completely controllable (resp. completely observable) systems are isomorphic to moduli spaces (in the sense of A. King \cite{King}) of $\theta$-stable representations of dimension vector $(1,n)$ for this quiver, where $\theta = (-n,1)$ (resp. $\theta = (n,-1)$). By computing the cohomology of these moduli spaces, as in \cite{Reineke}, we were then led to 

\begin{theorem} The moduli space $\wis{sys}_{m,n,p}^{cc}$ of completely controllable systems has a cell decomposition identical to the natural cell decomposition of a vectorbundle of rank $(p+1)n$ on the Grassmannian $\wis{Gras}_n(m+n-1)$ with respect to the Schubert cells on the Grassmannian.
\end{theorem}

In an earlier version of this note we claimed that the moduli space itself is a vectorbundle over the Grassmannian. However, this cannot be the case when $m=n$ as the referee kindly pointed out.

\subsection*{Acknowledgments:}

This paper was written while the second author enjoys a research stay at the University of Antwerp, with the aid of a grant of the European Science Foundation in the frame of the Priority Programme ``Noncommutative Geometry". We thank the anonymous referee for pointing out a mistake in the first version of this note.

\section{Proof of theorem 2}

Consider the quiver setting $(Q,\alpha)$ where the dimension vector is $\alpha = (1,n)$ and the quiver $Q$
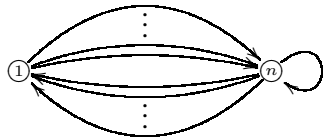
\[
\xymatrix@=3cm{
\vtx{1} \ar@/^/[r] \ar@/^2ex/[r]^{\vdots} \ar@/^5ex/[r] & \vtx{n} \ar@(ur,dr) \ar@/^/[l] \ar@/^2ex/[l] \ar@/^5ex/[l]_{\vdots} }
\]
has $m$ arrows $\{ b_1,\hdots,b_m \}$ from left to right and $p$ arrows $\{ c_1,\hdots,c_p \}$ from right to left. We can identify $V_{m,n,p}$ with $\wis{rep}_{\alpha}~Q$, where we associate to a system $\Sigma = (A,B,C)$ the representation $V_{\Sigma}$  which assigns to the arrow $b_i$ (resp. $c_j$) the $i$-th column $B_i$ of $B$ (resp. the $j$-th row $C^j$ of $C$) and the matrix $A$ to the loop. The basechange action of $(\lambda,g) \in GL(\alpha) = \C^* \times GL_n$ acts on the representation $V_{\Sigma} = (A,B_1,\hdots,B_m,C^1,\hdots,C^p)$ as follows:
\[
(\lambda.g).V_{\Sigma} = (gAg^{-1},gB_1 \lambda^{-1},\hdots,gB_m\lambda^{-1},\lambda C^1 g^{-1},\hdots,\lambda C^p g^{-1}), \]
and as the central subgroup $\C^*(1,\mathbf{1}_n)$ acts trivially on $\wis{rep}_{\alpha}~Q$, there is a natural one-to-one correspondence between equivalence classes of systems in $V_{m,n,p}$ and isomorphism classes of $\alpha$-dimensional representations in $\wis{rep}_{\alpha}~Q$. If $\C Q$ denotes the {\em path algebra} of the quiver $Q$, then it is well known that $\C Q$ is a formally smooth algebra in the sense of \cite{CuntzQuillen}, and that there is an equivalence of categories between finite dimensional right $\C Q$-modules and representations of $Q$. It is perhaps surprising that the system theoretic notion of canonical system corresponds under these identifications to the algebraic notion of simple module.

\begin{lemma} \label{lemma1} The following are equivalent:
\begin{enumerate}
\item{$\Sigma = (A,B,C) \in V_{m,n,p}$ is a canonical system,}
\item{$V_{\Sigma} = (A,B_1,\hdots,B_m,C^1,\hdots,C^p) \in \wis{rep}_{\alpha}~Q$ is a simple representation.}
\end{enumerate}
\end{lemma}

\Proof
$1 \Rightarrow 2$ : If $V_{\Sigma}$ has a proper subrepresentation of dimension vector $\beta = (1,l)$ for some $l < n$, then the rank of the control-matrix $c(\Sigma)$ is at most $l$, contradicting complete controllability. If $V_{\Sigma}$ has a proper subrepresentation of dimension vector $\beta' = (0,l)$ with $l \not= 0$, then the observation-matrix $o(\Sigma)$ has rank at most $n-l$, contradicting complete observability. $2 \Rightarrow 1$ : If $rk~c(\Sigma) = l < n$ then there is a proper subrepresentation of dimension vector $(1,l)$ of $V_{\Sigma}$. If $rk~o(\Sigma) = n-l$ with $l > 0$, then there is a proper subrepresentation of dimension vector $(0,l)$ of $V_{\Sigma}$.

\par \vskip 4mm
From \cite{LBProcesi} we recall that for a general quiver setting $(Q,\alpha)$ the isomorphism classes of $\alpha$-dimensional semi-simple representations are classified by the {\em affine} algebraic quotient variety
\[
\wis{rep}_{\alpha}~Q // GL(\alpha) = \wis{iss}_{\alpha}~Q \]
whose coordinate ring is generated by all traces along oriented cycles in the quiver $Q$. If $\alpha$ is the dimension vector of a simple representation, this affine quotient has dimension $1 - \chi_Q(\alpha,\alpha)$ where $\chi_Q$ is the {\em Euler form} of $Q$. Moreover, the isomorphism classes of {\em simple} representations form a Zariski open {\em smooth subvariety} of $\wis{iss}_{\alpha}~Q$. Specializing these general results from \cite{LBProcesi} to the case of interest, we recover Hazewinkels theorem.

\begin{theorem}[Hazewinkel] The moduli space $\wis{sys}_{m,n,p}^c$ of canonical systems is a smooth quasi-affine variety of dimension $(m+p)n$.
\end{theorem}

In fact, combining the theory of local quivers (see for example \cite{LBnagatn}) with the classification of all quiver settings having a smooth quotient variety due to Raf Bocklandt \cite{Bocklandt}, it follows that (unless $m=p=1$) $\wis{sys}_{m,n,p}^c$ is precisely the smooth locus of the affine quotient variety $\wis{iss}_{\alpha}~Q$.

\section{Proof of theorem 3}

For $(Q,\alpha)$  a quiver setting on $k$ vertices and if $\theta \in \Z^k$, a representation $V \in \wis{rep}_{\alpha}~Q$ is said to be {\em $\theta$-semistable} (resp. {\em $\theta$-stable}) if and only if for every {\em proper} non-zero subrepresentation $W$ of $V$ we have that $\theta.\beta \geq 0$ (resp. $\theta.\beta > 0$), where $\beta$ is the dimension vector of $W$. In the special case when $\alpha = (1,n)$ and $Q$ is the quiver introduced before, there are essentially two different {\em stability structures} on $\wis{rep}_{\alpha}~Q$ determined by the integral vectors
\[
\theta_+ = (-n,1) \qquad \text{and} \qquad \theta_- = (n,-1) \]
By the identification of $\wis{rep}_{\alpha}~Q$ with $V_{m,n,p}$ and the proof of lemma~\ref{lemma1} we have

\begin{lemma} \label{lemma2} For $\theta_+ = (-n,1)$ the following are equivalent:
\begin{enumerate}
\item{$\Sigma \in V_{m,n,p}$ is completely controllable,}
\item{$V_{\Sigma} \in \wis{rep}_{\alpha}~Q$ is $\theta_+$-stable.}
\end{enumerate}
For $\theta_- = (n,-1)$ the following are equivalent:
\begin{enumerate}
\item{$\Sigma \in V_{m,n,p}$ is completely controllable,}
\item{$V_{\Sigma} \in \wis{rep}_{\alpha}~Q$ is $\theta_-$-stable.}
\end{enumerate}
\end{lemma}

For a general stability structure $\theta$ and quiver setting $(Q,\alpha)$, A. King \cite{King} introduced and studied the {\em moduli space} $\wis{moduli}_{\alpha}^{\theta}~Q$ of $\theta$-semistable representations, the points of which classify isomorphism classes of direct sums of $\theta$-stable representations. In the case of interest to us we have
\[
\wis{sys}^{cc}_{m,n,p} = \wis{moduli}^{\theta_+}_{\alpha}~Q \qquad \text{and} \qquad
\wis{sys}^{co}_{m,n,p} = \wis{moduli}^{\theta_-}_{\alpha}~Q. \]
In \cite{Reineke} the Harder-Narasinham filtration associated to a stability structure was used to compute the cohomology of the moduli spaces $\wis{moduli}_{\alpha}^{\theta}~Q$ (at least if the quiver $Q$ has no oriented cycles). For general quivers the same methods can be applied to compute the number of $\FF_q$-points of these moduli spaces, where $\FF_q$ is the finite field of $q=p^l$ elements. In the case of interest to us, we get the rational functions
\[
\begin{cases}
\#~\wis{moduli}^{\theta_+}_{\alpha}~Q~(\FF_q) &= q^{n(p+1)} \prod_{i=1}^n \frac{q^{m+i-1}-1}{q^i-1} \\
\\
\#~\wis{moduli}^{\theta_-}_{\alpha}~Q~(\FF_q) &= q^{n(m+1)} \prod_{i=1}^n \frac{q^{p+i-1}-1}{q^i-1}, 
\end{cases}
\]
which suggests that the moduli space $\wis{sys}^{cc}_{m,n,p}$ is a vectorbundle of rank $n(p+1)$ over the Grassmannian $\wis{Gras}_n(m+n-1)$, and that the moduli space $\wis{sys}^{co}_{m,n,p}$ is a vectorbundle of rank $n(m+1)$ over $\wis{Gras}_n(p+n-1)$.

To a completely controllable $\Sigma = (A,B,C)$ one associates its {\em Kalman code} $K_{\Sigma}$, which is an array of $n \times m$ boxes $\{ (i,j)~|~0 \leq i < n,1\leq j\leq \}$, ordered lexicographically, with exactly $n$ boxes painted black. If the column $A^iB_j$ is linearly independent of all column vectors $A^kB_l$ with $(k,l) < (i,j)$ we paint box $(i,j)$ black. From this rule it is clear that if $(i,j)$ is a black box so are $(i',j)$ for all $i' \leq i$. That is, the Kalman code $K_{\Sigma}$ (which only depends on the $GL_n$-orbit of $\Sigma$) looks like
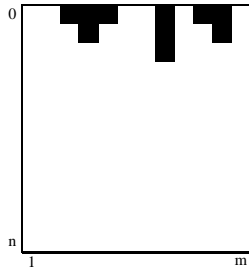
\[
\begin{xy} 0;/r.06pc/:
(60,65) *+(120,130){},{*\frm{**}};
(25,125) *+(10,10){},{*\frm{*}};
(35,120) *+(10,20){},{*\frm{*}};
(45,125) *+(10,10){},{*\frm{*}};
(75,115) *+(10,30){},{*\frm{*}};
(95,125) *+(10,10){},{*\frm{*}};
(105,120) *+(10,20){},{*\frm{*}};
(-5,125) *+{\txt{\tiny 0}};
(-5,5) *+{\txt{\tiny n}};
(5,-5) *+{\txt{\tiny 1}};
(115,-5) *+{\txt{\tiny m}};
\end{xy}
\]
Assume $\kappa = K_{\Sigma}$ has $k$ black boxes on its first row at places $(0,i_1),\hdots,(0,i_k)$. Then we assign to $\kappa$ the strictly increasing sequence 
\[
1 \leq j_{\kappa}(1)=i_1 < j_{\kappa}(2)=i_2 < \hdots < j_{\kappa}(k)=i_k \leq m \]
and another sequence $p_{\kappa}(1),\hdots,p_{\kappa}(k)$, where $p_{\kappa}(j)$ is the total number of black boxes in the $i_j$-th column of $\kappa$, that is,
\[
p_{\kappa}(1) + p_{\kappa}(2) + \hdots + p_{\kappa}(k) = n. \]
It is clear that there is a one-to-one correspondence between Kalman codes and pairs of functions satisfying these conditions. Further, define the strictly increasing sequence
\[
h_{\kappa}(0)=0 <  h_{\kappa}(1)=p_{\kappa}(1) < \hdots <  h_{\kappa}(j) = \sum_{i=1}^j p_{\kappa}(i) < \hdots < h_{\kappa}(k) = n. \]
With these notations we have the following canonical form for $\Sigma = (A,B,C) \in V_{m,n,p}^{cc}$ which is essentially \cite[lemma 3.2]{Geiss}:

\begin{lemma} For a completely reachable system $\Sigma = (A,B,C)$ with Kalman code $\kappa = K_{\Sigma}$, there is a unique $g \in GL_n$ such that
$g.(A,B,C) = (A',B',C')$
with
\begin{itemize}
\item{$B'_{j_{\kappa}(i)} = \mathbf{1}_{h_{\kappa}(i-1)+1}$ for all $1 \leq i \leq k$.}
\item{$A'_i = \mathbf{1}_{i+1}$ for all $i \notin \{ h_{\kappa}(1),h_{\kappa}(2),\hdots,h_{\kappa}(k) \}$.}
\item{All entries in the remaining columns of $A'$ and $B'$ are determined as the quotient of two specific $n \times n$ minors of $c(\Sigma)$.}
\item{$C' = Cg^{-1}$.}
\end{itemize}
\end{lemma}

This allows us to prove theorem 3 :

\begin{theorem} The moduli space $\wis{sys}_{m,n,p}^{cc}$  of completely controllable systems has a cell decomposition identical to the natural cell decomposition of a vectorbundle of rank $n(p+1)$ over the Grassmann manifold $\wis{Gras}_n(m+n-1)$.
\end{theorem}

\Proof
Define a map $V_{m,n,p}^{cc} \rTo^{\phi} \wis{Gras}_n(m+n-1)$ by sending a completely reachable system $\Sigma = (A,B,C)$ to the point in $\wis{Gras}_n(m+n-1)$ determined by the $n \times (m+n-1)$ matrix
\[
M_{\Sigma} = \begin{bmatrix} B'_1&  \hdots  & B'_m & A'_1 &  \hdots  & A'_{n-1} \end{bmatrix},
\]
where $(A',B',C')$ is the canonical form of $\Sigma$ given by the previous lemma. By construction,
$M_{\Sigma}$ has rank $n$ with invertible $n \times n$ matrix determined by the columns
\[
I_{\kappa} = \{ j_{\kappa}(1) < \hdots < j_{\kappa}(k) < m+c_1 < \hdots < m+c_{n-k} \} \subset \{ 1,\hdots,m+n-1 \}, \]
where $\{ c_1,\hdots,c_{n-k} \} = \{ 1,\hdots,n \} - \{ h_{\kappa}(1),\hdots,h_{\kappa}(k) \}$. As all remaining entries of $(A',B')$ are determined by $c(\Sigma)$ it follows that $\phi(\Sigma)$ depends only on the $GL_n$-orbit of $\Sigma$, whence the map factorizes through
\[
\wis{sys}_{m,n,p}^{cc} \rTo^{\psi} \wis{Gras}_n(m+n-1), \]
and we claim that $\psi$ is surjective. To begin, all multi-indices $I=\{ 1 \leq  d_1 < d_2 < \hdots < d_n \leq n+m-1 \}$ are of the form $I_{\kappa}$ for some Kalman code $\kappa$. Define
\[
\{ d_1,\hdots,d_n \} = \{ i_1,\hdots,i_k \} \cup \{ m+c_1,\hdots,m+c_{n-k} \} \]
with $i_j \leq m$ and $1 \leq c_j < n$, and let $\{ e_1 < \hdots < e_k \} = \{ 1,\hdots,n \} - \{ c_1,\hdots,c_{n-k} \}$, and set $e_{0}=0$. Construct the Kalman code $\kappa$ having $e_{j}-e_{j-1}$ black boxes in the $i_j$-th column and verify that $I$ is indeed $I_{\kappa}$.

$\wis{Gras}_n(m+n-1)$ is covered by modified Schubert cells $S_I$ (isomorphic to some affine space) consisting of points such that the $I$-minor is invertible, where $I$ is a multi-index $\{ d_1,\hdots,d_n \}$, and the dimension of the subspace spanned by the first $k$ columns is $i$ iff $k < d_{i+1}$. A point in $S_I$ can be taken such that the $d_i$-th column is equal to
\[
\begin{cases}
\mathbf{1}_{h_{\kappa}(i-1)+1} & \qquad  \txt{for $d_i \leq m$} \\
\mathbf{1}_{j+1} & \qquad \txt{for $d_i = m+i$},
\end{cases}
\]
where $I = I_{\kappa}$. This determines a $n \times (n+m-1)$ matrix $\begin{bmatrix} B_1 & \hdots & B_m & A_1 & \hdots & A_{n-1} \end{bmatrix}$, and choosing any last column $A_n$ and any $p \times n$ matrix $C$ we obtain a system
$\Sigma = (A,B,C)$ which is completely controllable, and which is mapped to the given point under $\psi$. This finishes the proof.

\par \vskip 4mm

Because the map $(A,B,C) \rTo (A^{tr},C^{tr},B^{tr})$ defines a duality between $V_{m,n,p}^{co}$ and $V_{p,n,m}^{cc}$, we have a similar result for the moduli spaces of completely observable systems.

\begin{theorem} The moduli space of completely observable systems $\wis{sys}^{co}_{m,n,p}$ has a cell decomposition identical to that of a vectorbundle of rank $n(p+1)$ over the Grassmann manifold $\wis{Gras}_n(p+n-1)$.
\end{theorem}

\section{Proof of theorem 1}

The counting argument of the previous section gives us also a conjectural description of the infinite dimensional variety admitting a stratification by the moduli spaces $\wis{sys}_{m,n,p}^{cc}$. It follows from the explicit rational form of $\#~\wis{sys}_{m,n,p}^{cc}~(\FF_q)$ and the {\em $q$-binomial theorem} that 
\[
\sum_{n=0}^{\infty}~\#~\wis{sys}_{m,n,p}^{cc}~(\FF_q)~t^n = \prod_{i=1}^m \frac{1}{1-q^{p+i}t} \]
In the special case when $p=0$ we recover the cohomology of the infinite Grassmannian $\wis{Gras}_m(\infty)$ of $m$-dimensional subspaces of a countably infinite dimensional vectorspace. For $p \geq 1$ we only get a factor of the cohomology of $\wis{Gras}_{m+p}(\infty)$, which led to the following result.

\begin{theorem} The disjoint union $\bigsqcup_n~\wis{sys}_{m,n,p}^{cc}$ is the open subset of the
infinite dimensional Grassmann manifold $\wis{Gras}_{m+p}(\infty)$ which is the union of all standard affine open sets corresponding to a multi-index set $I = \{ 1 \leq d_1 < d_2 < \hdots < d_{m+p} \}$ such that
\[
\{ m+1,m+2,\hdots,m+p,m+p+n \} \subset I. \]
\end{theorem}

\Proof
Let $\Sigma = (A,B,C)$ be a completely controllable system in canonical form represented by the point $p_{\Sigma} \in \wis{sys}_{m,n,p}^{cc}$. Consider the $n \times (m+p+n)$ matrix
\[
L_{\Sigma} = \begin{bmatrix} B & C^{tr} & A \end{bmatrix}. \]
The submatrix $M_{\Sigma} = \begin{bmatrix} B_1 & \hdots & B_m & A_1 & \hdots & A_{n-1} \end{bmatrix}$ has rank $n$, whence so has $L_{\Sigma}$, and $p_{\Sigma}$ determines a point in $\wis{Gras}_n(n+m+p)$. Under the natural duality
\[
\wis{Gras}_n(m+p+n) \rTo^D \wis{Gras}_{m+p}(m+p+n), \]
the point $p_{\Sigma}$ is mapped to the point determined by the $(m+p) \times (m+p+n)$ matrix $N_{\Sigma}$ whose rows give a basis for the linear relations holding among the columns of $L_{\Sigma}$. Because $M_{\Sigma}$ has rank $n$ it follows that the columns of $C^{tr}$ and the last column $A_n$ of $A$ are linearly dependent of those of $M_{\Sigma}$. As a consequence the matrix
\[
N_{\Sigma} = \begin{bmatrix}
U_1 & \hdots & U_m & V_1 & \hdots & V_p & W_1 & \hdots & W_n \end{bmatrix} \]
has the property that the submatrix $\begin{bmatrix} V_1 & \hdots & V_p & W_n \end{bmatrix}$ has rank $p+1$. This procedure defines a morphism
\[
\wis{sys}_{m,n,p}^{cc} \rTo^{\gamma_n} \wis{Gras}_{m+p}(m+p+n), \]
the image of which is the open union of all standard affine opens determined by a multi-index set $I = \{ 1 \leq d_1 < d_2 < \hdots < d_{m+p} \leq m+p+n \}$ satisfying
\[
\{ m+1,m+2,\hdots,m+p,m+p+n \} \subset I. \]
Therefore, the image of the morphism
\[
\bigsqcup_n~\wis{sys}_{m,n,p}^{cc} \rTo^{\sqcup \gamma_n} \wis{Gras}_{m+p}(\infty) \]
is the one of the statement of the theorem. The dimension $n$ of the system corresponding to a point in this open set of $\wis{Gras}_{m+p}(\infty)$ is determined by $d_{m+p} = m+p+n$.

\par \vskip 4mm

By the duality between $V_{m,n,p}^{cc}$ and $V_{p,n,m}^{co}$ used in the previous section we deduce:

\begin{theorem} The disjoint union $\bigsqcup_n~\wis{sys}_{m,n,p}^{co}$ is the open subset of $\wis{Gras}_{m+p}(\infty)$ which is the union of all standard affine opens 
corresponding to a multi-index set $I = \{ 1 \leq d_1 < d_2 < \hdots < d_{m+p} \}$ such that
\[
\{ 1,2,\hdots,m,m+p+n \} \subset I. \]
\end{theorem}

This, in turn, proves theorem~1 :

\begin{theorem} The disjoint union $\bigsqcup_n~\wis{sys}^c_{m,n,p}$ of all moduli spaces of canonical systems with fixed input- and output-dimension $m$ and $p$ is the open subset of the infinite Grassmannian $\wis{Gras}_{m+p}(\infty)$ of $m+p$-dimensional subspaces of a countably infinite dimensional vectorspace which is the intersection of all possible standard open subsets $X_I$ and $X_J$, where $I$ and $J$ are multi-index sets satisfying the conditions
\[
\{ m+1,m+2,\hdots,m+p,m+p+n \} \subset I \qquad \text{and} \qquad \{ 1,2,\hdots,m,m+p+n \} \subset J. \]
\end{theorem}

\providecommand{\href}[2]{#2}
\begingroup\raggedright\endgroup

\end{document}